\RequirePackage[2020-02-02]{latexrelease}
\documentclass[11pt]{amsart}
\usepackage{amsfonts}
\usepackage{amsmath}
\usepackage{amsthm}
\usepackage{amssymb}
\usepackage{wasysym}
\usepackage[numbers]{natbib}
\usepackage[fit]{truncate}
\usepackage{hyperref}
\usepackage{theoremref}

\usepackage{mathtools}

\usepackage[dvipsnames]{xcolor}

\usepackage{tikz-cd}

\usepackage{thmtools}
\usepackage{thm-restate}

\usepackage{fmtcount}

\usepackage{enumitem}

\usepackage[UKenglish]{babel}
\usepackage[UKenglish]{isodate}

\usepackage{verbatim}
\usepackage[utf8]{inputenc}

\usepackage{dsfont}

\numberwithin{equation}{section}

\newenvironment{acknowledgements}{\ 
	
	{\textsl{Acknowledgements.}}}{}

\makeatletter
\def\blfootnote{\gdef\@thefnmark{}\@footnotetext}
\makeatother

\newcommand{\N}[0]{\mathbb{N}}
\newcommand{\Z}[0]{\mathbb{Z}}
\newcommand{\Q}[0]{\mathbb{Q}}
\newcommand{\R}[0]{\mathbb{R}}

\newcommand{\MM}[0]{\mathcal{M}}

\renewcommand{\AA}[0]{\mathcal{A}}
\newcommand{\BB}[0]{\mathcal{B}}
\newcommand{\KK}[0]{\mathcal{K}}

\newcommand{\brackets}[1]{\left( #1 \right)}

\newcommand{\LL}[0]{\mathcal{L}}
\newcommand{\Lor}[0]{\mathcal{L}_{\mathrm{or}}}

\newcommand{\Lexp}[0]{\mathcal{L}_{\exp}}
\newcommand{\Trcf}[0]{T_{\mathrm{rcf}}}

\newcommand{\Texp}[0]{T_{\exp}}

\newcommand{\ul}[1]{\underline{#1}}

\newcommand{\Kexp}[0]{(K,E)}
\newcommand{\Rexp}[0]{\mathbb{R}_{\exp}}

\DeclareMathOperator{\dcl}{dcl}
\DeclareMathOperator{\Th}{Th}

\DeclareMathOperator{\EXP}{EXP}
\DeclareMathOperator{\ecl}{CL}

\theoremstyle{definition}
\newtheorem{defn}{Definition}[section]
\newtheorem{definition}[defn]{Definition}

\theoremstyle{plain}
\newtheorem{lemma}[defn]{Lem\-ma}
\newtheorem{observation}[defn]{Observation}
\theoremstyle{plain}
\newtheorem*{lemma*}{Lem\-ma}
\theoremstyle{plain}
\newtheorem{proposition}[defn]{Proposition}
\newtheorem{fact}[defn]{Fact}

\newtheorem{question}[defn]{Question}
\newtheorem{corollary}[defn]{Corollary}
\theoremstyle{plain}
\theoremstyle{plain}
\newtheorem{theorem}[defn]{Theorem}
\newtheorem*{schanuel}{Schanuel's Conjecture}
\newtheorem*{transfer}{Transfer Conjecture}
\newtheorem*{konig}{K\H{o}nig's Lemma}
\theoremstyle{definition}

\theoremstyle{definition}

\theoremstyle{plain}

\begin{document}
	
	\title[Embedding the prime model of real exponentiation]{Embedding the prime model of real exponentiation into o-minimal exponential fields}
	
	\author[L.~S.~Krapp]{Lothar Sebastian Krapp}
	
	\address{Fachbereich Mathematik und Statistik, Universität Konstanz, 78457 Konstanz, Germany}
	\email{sebastian.krapp@uni-konstanz.de}
	
	\date{\today}
	
	\cleanlookdateon
	
	\maketitle
	
	\begin{abstract}
		Motivated by the decidability question for the theory of real exponentiation and by the Transfer Conjecture for o-minimal exponential fields, 
		we show that, under the assumption of Scha\-nu\-el's Conjecture, the prime model of real exponentiation is embeddable into any o-minimal exponential field, {\color{black} where the embedding is not necessarily elementary.
		This is a consequence of an unconditional 
		model theoretic embeddability result that we obtain by applying K\H{o}nig's Lemma.}
	\end{abstract}
	
	
	\blfootnote{\textup{2020} \textit{Mathematics Subject Classification}: 03C64 (03C60 12J15 12L12)
		}
	
	\section{Introduction}
	
	In his highly influential work \cite{tarski}, Tarski proved that the complete theory $\Trcf$ of the real closed field $\R$ in the language of ordered rings $\Lor=\{+,-,\cdot,0,1,<\}$ is decidable, by presenting an explicit quantifier-elimination algorithm for $\Trcf$. As a result, he noted that the $\Lor$-structure $(\R,+,-,\cdot,0,\allowbreak1,\allowbreak<)$ is -- in modern terminology -- o-minimal, i.e.\ any unary definable subset of this structure is a finite union of points and open intervals. In the same work, Tarski asked whether decidability can also be obtained for the complete theory $\Texp$ of the real exponential field $\Rexp = (\R,+,-,\cdot,0,1,\allowbreak <,\exp)$, where $\exp$ denotes the standard exponential function $x\mapsto \mathrm{e}^x$ on $\R$. While this question is open to the date, considerable progress has been made since the 1990s: In \cite{wilkie}, Wilkie proved that $\Rexp$ is model complete, and thus this structure is o-minimal (see also van den Dries, Macintyre and Marker~\cite{dries2}). Bulding on this result, Macintyre and Wilkie \cite{macintyre} showed that under the assumption of the real version of Schanuel's Conjecture below, $\Texp$ is decidable.
	\begin{schanuel}\thlabel{schnauel}
		Let $n\in\N$ and let $\alpha_1,\ldots,\alpha_n\in \R$ be $\Q$-linearly independent. Then the transcendence degree of 		$\Q(\alpha_1,\ldots,\alpha_n,\mathrm{e}^{\alpha_1},\ldots,\mathrm{e}^{\alpha_n})$ over $\Q$ is at least $n$.
	\end{schanuel}
	In general, an ordered exponential field $\KK=(K,+,-,\cdot,0,1,<,E)$ is an expansion of an ordered field $(K,+,-,\cdot,0,1,<)$ by an exponential $E$, i.e.\ an order-preserving isomorphism from the ordered additive group $(K,+,0,<)$ to the ordered multiplicative group $(K^{>0},\cdot,1,<)$. We denote the corresponding first-order language $\Lor\cup \{E\}$, where $E$ is a unary function symbol, by $\Lexp$. Following the terminology of Krapp~\cite{krappla}, we call an ordered exponential field $\KK$ an $\EXP$-field if its exponential satisfies the first-order $\Lexp$-sentence expressing the differential equation $E'=E$ {\color{black} with initial condition $E(0)=1$}. While $\Rexp$ and, more generally, any model of $\Texp$ is an o-minimal $\EXP$-field, {\color{black} the following is still open:}
	{\color{black}}
	\begin{transfer}\thlabel{transfer}
		Any o-minimal $\EXP$-field is elementarily equivalent to $\Rexp$.
	\end{transfer}
	Berarducci and Servi~\cite{berarducci} showed that the Transfer Conjecture would imply that $\Texp$ is decidable, {\color{black}motivating the study of o-minimal $\EXP$-fields.} 
	One approach towards proving the Transfer Conjecture is to show that the (unique) prime model of real exponentiation, i.e.\ the prime model of $\Texp$, elementarily embeds into any o-minimal $\EXP$-field, which is also the motivating question for this line of research.
	
	\begin{question}[Main question]\thlabel{qu:main}
		Does the prime model of $\Texp$ elementarily embed into any o-minimal $\EXP$-field?
	\end{question}
	
	{\color{black} Since} \Autoref{qu:main} can be answered positively if and only if the Transfer Conjecture holds, 
	it is natural to {\color{black} approach this question} only under the assumption of Schanuel's Conjecture. 
	While this note does not provide a complete answer to \Autoref{qu:main}, we show in \Autoref{cor:scprimeemb} that under the assumption of Schanuel's Conjecture, the prime model of $\Texp$ embeds into any o-minimal $\EXP$-field (however, this embedding might not necessarily be elementary). 
	{\color{black}\Autoref{cor:scprimeemb} is deduced from the model theoretic main result of this paper (\Autoref{thm:eclembeddable1}) which establishes sufficient conditions on two structures $\AA$ and $\BB$ and a substructure $\AA'$ of $\AA$ in order that $\AA'$ embeds into $\BB$. By application of \Autoref{thm:eclembeddable1}, we also obtain an} embeddability result (\Autoref{thm:eclembeddable}) for 
	exponential algebraic closures of $\Z$ within one o-minimal $\EXP$-field into another, {\color{black} without assuming Schanuel's Conjecture}.

	\subsection{Notation and terminology}
		More background on the model theoretic notation and terminology we use can be found in Marker~\cite{marker}. If it is clear from the context, then the $\Lor$-structure of an ordered field $(K,+,-,\cdot,0,1,<)$ is simply denoted by $K$ and the $\Lexp$-structure of an ordered exponential field $(K,+,-,\cdot,0,1,<,E)$ simply by $(K,E)$.
		Given an ordered exponential field $(K,E)$, we say that a subfield $F\subseteq K$ is \textbf{exponentially closed} in $(K,E)$ if the restriction $E|_{F}$ is an exponential on $F$. In this case we also denote $E|_F$ simply by $E$. 
		For any structure $\mathcal{M}=(M,\ldots)$, its complete theory is denoted by $\Th(\MM)$ and its existential theory by $\Th_\exists(\MM)$. For instance, $\Texp = \Th(\Rexp)$ denotes the theory of real exponentiation and $\Th_\exists(\Rexp)$ denotes the existential theory of real exponentiation. 
		The definable closure of a set $C$ in $\MM$ is denoted by $\dcl(C;\MM)$. We say that a set is definable if it is definable \emph{with} parameters. 
		{\color{black} For a formula $\varphi(x_1,\ldots,x_n)$ we denote the subset of $M^n$ it defines by $\varphi(\MM)=\{(a_1,\ldots,a_n)\in M^n\mid \MM\models \varphi(a_1,\ldots,a_n)\}$.}
		Variable tuples are denoted by $\ul{x}$, and we only specify their length if it is of importance. 
		If $\ul{x}=(x_1,\ldots,x_n)$ and $f$ is a unary map, we write $f(\ul{x})$ for the tuple $(f(x_1),\ldots,f(x_n))$. 
		We denote by $\N$ the set of natural numbers \emph{without} $0$.
		\textbf{Throughout the rest of this note, let $(K,E)$ denote an o-minimal $\EXP$-field.}

\section{Embedding substructures}

	We start by proving our model theoretic main result relying on K\H{o}nig's Lemma, which we briefly recall in the following (see Marker~\cite[\linebreak page~320~f.]{marker}).
	A \textbf{finite branching tree} is a partial order $(T,<)$ that contains some $r\in T$ such that for any $t\in T$ we have $r\leq t$, $\{s\in T\mid s<t\}$ is linearly ordered by $<$ and $t$ only has finitely many immediate successors in $T$ (i.e.\ finitely many $s>t$ such that $T$ contains no element strictly between $t$ and $s$). A \textbf{path} through $T$ is a function $f\colon \omega\to T$ such that $f(n)<f(n+1)$ for any $n<\omega$.
	
	\begin{konig}
		Let $(T,<)$ be an infinite finite branching tree. Then there is a path through $T$.
	\end{konig}
	
	\begin{theorem}\thlabel{thm:eclembeddable1}
		Let $\LL$ be a countable language, let $\AA$ and $\BB$ be two $\LL$-structures and let $\AA'$ be a countably infinite substructure of $\AA$. Suppose that the following hold:
		\begin{enumerate}[label = (\roman*)]
			\item $\BB\models \Th_\exists(\AA)$.
			
			\item For any $a\in A'$, there is an existential $\LL$-formula $\theta(x)$ such that $\AA\models \theta(a)$ and both $\theta(\AA)$ and $\theta(\BB)$ are finite (not necessarily of equal cardinality).
		\end{enumerate}
		Then $\AA'$ can be embedded into $\BB$ as an $\LL$-substructure.
	\end{theorem}

	\begin{proof}
		Let $(a_n)_{n<\omega}$ be an enumeration of $A'$ and, for any $n<\omega$, let $\theta_n(x)$ be an existential $\LL$-formula such that $\AA\models \theta_n(a_n)$ and $|\theta_n(\AA)|,|\theta_n(\BB)|<\infty$. Note that $\theta_n(\BB)\neq\emptyset$, as $\BB\models \Th_\exists(\AA)$.
		Let $(\exists \ul{x}\ \psi_n(\ul{x}))_{n<\omega}$ be an enumeration of $\Th_\exists(\AA')$, where each $\psi_n(\ul{x})$ is quantifier-free. For our later argument, we may assume that the free variables of each $\psi_n(\ul{x})$ are contained in $\{x_0,\ldots,x_n\}$ and that $\AA\models \psi_n(a_0,\ldots,a_n)$, as otherwise we may reorder the enumeration and relabel the free variables in each formula.
				
		For each $n<\omega$, let $\psi_n'(\ul{x})$ be the conjunction of $\psi_n(\ul{x})$ with $\theta_0(x_0)\wedge\ldots \wedge \theta_n(x_n)$. 
		Then $\AA\models \psi_n'(\ul{a})$ and $\BB\models \exists \ul{x}\ \psi_n'(\ul{x})$, as $\exists \ul{x}\ \psi_n'(\ul{x})$ is logically equivalent to a sentence in $\Th_\exists(\AA)$. Thus, both $\psi_n'(\AA)$ and $\psi_n'(\BB)$ are non-empty and finite.
		Moreover, let $T_n\subseteq B^{n+1}$ consist of all tuples $\ul{b}=(b_0,\ldots,b_n)\in B^{n+1}$ with $\BB\models \psi_0'(\ul{b})\wedge\ldots \wedge \psi_n'(\ul{b})$. Again, $T_n$ is non-empty, as $\BB\models \Th_\exists(\AA)$, as well as finite. 
		
		Let $T$ be the disjoint union $T=\bigcup_{n=-1}^\infty T_n$, where $T_{-1}:=\{r\}$ with $r=(\ )$ denoting the null tuple. For any $i,j\in \{-1,0,1,\ldots\}$ with $i\leq j$ and any $s\in T_i$ and $t\in T_j$ we set $s\leq t$ if $s$ is the projection of $t$ onto its first $i+1$ components.
		Then $(T,<)$ is a finite branching tree with minimal element $r$. By K\H{o}nig's Lemma, there is some path $f\colon \omega\to T$ through $T$. For any $n<\omega$, we let $b_n$ be the $(n+1)$-th entry of the tuple $f(n+1)$. Then $\BB\models \psi_0'(\ul{b})\wedge\ldots \wedge \psi_n'(\ul{b})$.
		
		Finally, we let $\iota\colon A'\to B$ map each $a_n$ to $b_n$. Note that for any $n<\omega$ we have $\BB\models \psi_n(\ul{b})$. 
		Due to the enumeration we chose for $\Th_\exists(\AA')$, we thus obtain that $\iota$ preserves existential formulas, i.e.\ for any existential $\LL$-formula $\varphi(\ul{x})$ with $\AA'\models \varphi(\ul{a})$ also $\BB\models \varphi(\ul{b})$. Hence, $\iota$ preserves quantifier-free formulas in both directions, showing that $\iota$ is, indeed, an $\LL$-embedding.		
	\end{proof}

	\section{Exponential algebraic closures}
	
	In the following, we mostly follow the terminology of Macintyre~\cite{macintyre2} and Kirby~\cite{kirby} adjusted to our context. 
	For a subring $A\subseteq K$, we 
	say that $A$ is an \textbf{$E$-ring} if it is closed under $E$, i.e.\ for any $a\in A$ also $E(a)\in A$. 
	The smallest $E$-ring in $K$ is denoted by $\Z^{E}$. 
	If $B$ is a ring, then we	
	denote by $B[\ul{x}]_E=B[x_1,\ldots,x_n]_{E}$ the ring of \textbf{exponential polynomials} in $n$ variables, which consists of all functions of the form $p(\ul{x},E(\ul{x}))$ for some $p(\ul{x},\ul{y})\in B[x_1,\ldots,x_n,y_1,\ldots,y_n]$.\footnote{Note that \cite{macintyre2} and \cite{kirby} use free $E$-rings in several variables, which allow multiple iterations of exponentiation. However, for the purpose of this note it suffices to work with the ring of exponential polynomials (cf.\ Servi~\cite[Proposition~4.5.4]{servi}).}

	\begin{definition}
		\begin{enumerate}[label = (\roman*)]
			\item Let $n \in \N$ and let $A \subseteq K$ be an $E$-ring. 
			Moreover, let $f_1,\ldots,f_n \in A[x_1,\ldots,x_n]_{E}$ and let  $|J_{f_1,\ldots,f_n}(\ul{x})|$ denote the determinant of the Jacobian matrix $$J_{f_1,\ldots,f_n}(\ul{x}) = \brackets{\begin{array}{ccc}
					\frac{\partial f_1}{\partial x_1}(\ul{x})& \ldots & \frac{\partial f_1}{\partial x_n}(\ul{x})\\
					\vdots & & \vdots\\
					\frac{\partial f_n}{\partial x_1}(\ul{x})& \ldots & \frac{\partial f_n}{\partial x_n}(\ul{x})
			\end{array}}\!\!.$$
			Then the \textbf{Khovanskii system $S(\ul{x})$} of $f_1,\ldots,f_n$ over $A$ is the system of equations and inequations 
			$$f_1(\ul{x}) = \ldots = f_n(\ul{x}) = 0 \text{ and }|J_{f_1,\ldots,f_n}(\ul{x})| \neq 0.$$
		
			\item Let $B$ be a subset of $K$ and let $A\subseteq K$ be the smallest $E$-ring containing $B$. 
			An element $a_1 \in K$ is said to be \textbf{exponentially algebraic} over $B$ if for some $n \in \N$ there exist a Khovanskii system $S(x_1,\ldots,x_n)$ over $A$ and $a_2,\ldots,a_n \in K$ such that $(a_1,\ldots,a_n)$ solves this system. The set of all elements of $K$ that are exponentially algebraic over $B$ is called the \textbf{exponential algebraic closure} of $B$ in $(K,E)$ and is denoted by $\ecl^{E}_K(B)$.
		\end{enumerate}
	\end{definition}

	The following can be obtained by simple Jacobian calculations (see also Kirby~\cite[Lemma~3.3]{kirby} and Macintyre~\cite[Lemma~21]{macintyre2}).
	
	\begin{lemma}\thlabel{lem:expclosed}
		Let $B$ be a subset of $K$. Then $\ecl^{E}_K(B)$ is a subfield of $K$ that is exponentially closed in $(K,E)$.
	\end{lemma}
	
	Note that the condition on the determinant of the Jacobian to be non-zero implies that Khovanskii systems only have isolated solutions. This is a consequence of the Implicit Function Theorem in o-minimal structures (see van den Dries~\cite[page~113]{dries}).
	
	\begin{lemma}\thlabel{lem:khovanskiifinite}
		Let $A\subseteq K$ be an $E$-ring and let $S(\ul{x})$ be a Khovanskii system over $A$. Then $S$ only has finitely many solutions in $K$.
	\end{lemma}
	
	\begin{proof}
		It only remains to note that $S(\ul{x})= S(x_1,\ldots,x_n)$ defines finitely many connected components in $K^n$ (see \cite[\S~3.2]{dries}).
	\end{proof}
		
	We now investigate the connection between the definable closure and the exponential algebraic closure in o-minimal $\EXP$-fields. First, we make an observation following from \Autoref{lem:khovanskiifinite}.

	\begin{observation}\thlabel{prop:dclexpecl2}
		Let $B$ be a subset of $K$. Then $\ecl^{E}_K(B) \subseteq \dcl(B;\allowbreak\Kexp)$.
	\end{observation}

	Recall that due to definable Skolem functions in o-minimal expansions of ordered groups, for any subset $B$ of $K$ we have that  $(\dcl(B;\Kexp),E)\preceq (K,E)$, i.e.\ $\dcl(B;\Kexp)$ is the domain of an elementary substructure of $(K,E)$. Hence, $(\dcl(\emptyset;\Kexp),E)$ is the unique prime model of $\Th(K,E)$ (see Krapp~\cite[Proposition~4.75]{krappdiss} for further details).
	If $(K,E)$ is not only assumed to be an o-minimal $\EXP$-field but already a model of real exponentiation, then we can strengthen the conclusion of \Autoref{prop:dclexpecl2} to $\ecl^{E}_K(B) = \dcl(B;\allowbreak\Kexp)$. 
	This result is mentioned by Macintyre in \cite[Theorem~22]{macintyre2}, who attributes it to Wilkie~\cite{wilkie}.
	We point out that it is also an immediate consequence of Jones and Wilkie \cite[Theorem~4.2]{jones} applied to the locally polynomially bounded structure $\mathcal{M}=(K,\mathcal{F})$ with $\mathcal{F}=\{E\}$.
	
	\begin{proposition}\thlabel{prop:dclexpecl}\thlabel{prop:dclexpecl3}\thlabel{prop:primekhovanskii}
		Let $(K,E) \models \Texp$ and let $B$ be a subset of $K$. Then $$\dcl(B;(K, E)) = \ecl_K^{E}(B).$$ 
		In particular, $(\ecl_K^{E}(B),E)\preceq (K,E)$ and  $(\ecl_K^{E}(\emptyset),E)$ is the prime model of $\Texp$.
	\end{proposition}
	
	The aim of this note is to show that, under the assumption of Schanuel's Conjecture, the prime model of $\Texp$ is embeddable into any o-minimal $\EXP$-field (see \Autoref{cor:scprimeemb}). This result will be deduced from the following.
	
	\begin{theorem}\thlabel{thm:eclembeddable}
		Let $(K_1,E_1)$ and $(K_2,E_2)$ be o-minimal $\EXP$-fields. Suppose that $(K_2,E_2) \models \Th_\exists(K_1,E_1)$. Then there exists an $\Lexp$-em\-bed\-ding of the ordered exponential field $(\ecl^{E_1}_{K_1}(\emptyset),E_1)$ into $(K_2,E_2)$.
	\end{theorem}
	
	\begin{proof}
		We apply \Autoref{thm:eclembeddable1} to $\AA= (K_1,E_1)$, $\BB=(K_2,E_2)$ and $\AA'=(\ecl^{E_1}_{K_1}(\emptyset),E_1)$, which is countably infinite due to \Autoref{lem:khovanskiifinite}. 
		In order to do so, we let $a\in \ecl^{E_1}_{K_1}(\emptyset)$. Then there is a Khovanskii system $S(x,\ul{y})$ over $\Z^E$ such that $\AA\models \exists \ul{y}\ S(a,\ul{y})$. Setting $\theta(x)$ to be $ \exists \ul{y}\ S(x,\ul{y})$ we are done by \Autoref{lem:khovanskiifinite}, as $\mathcal{B}\models \Th_{\exists}(\mathcal{A})$. 
	\end{proof}

	\section{Prime model of real exponentiation}
	
	\textbf{In this section, we denote by $(P,\exp)$ the prime model of $\Texp$.}
	Since $(P,\exp)\models \Texp$, any structure containing $(P,\exp)$ as a substructure already satisfies the existential theory $\Th_\exists(\Rexp)$.
	Thus, by \Autoref{prop:primekhovanskii} and \Autoref{thm:eclembeddable} we obtain the following.
		
	\begin{corollary}\thlabel{cor:primemodelemb} 
		There is an emebedding of $(P,\exp)$ into $(K,E)$ if and only if $(K,E) \models \Th_\exists(\Rexp)$.
	\end{corollary}
	
	\Autoref{cor:primemodelemb} shows that o-minimal $\EXP$-fields satisfying the existential theory of real exponentiation contain the prime model of real exponentiation as a substructure. Under the assumption of Schanuel's Conjecture, any o-minimal $\EXP$-field satisfies $\Th_{\exists}(\Rexp)$. This fact is basically due to Servi~\cite[page~104]{servi} (see also Krapp~\cite[Proposition~4.57]{krappdiss}).

	\begin{fact}\thlabel{fact:servi}
		Assume Schanuel's Conjecture. Then any o-minimal $\EXP$-field satisfies $\Th_\exists(\Rexp)$.
	\end{fact}
	
	 As a result of \Autoref{cor:primemodelemb} and \Autoref{fact:servi}, we obtain the desired main result of this note.
	
	\begin{theorem}\thlabel{cor:scprimeemb}
		Assume Schanuel's Conjecture. Then $(P,\exp)$ embeds into any o-minimal $\EXP$-field.
	\end{theorem}

\begin{acknowledgements}
	\Autoref{thm:eclembeddable} and \Autoref{cor:scprimeemb} were part of my doctoral thesis \cite{krappdiss}. In this regard, I wish to thank my supervisor Salma Kuhlmann for her continuous help and support. I also thank Alessandro Berarducci for an insightful discussion about \cite[Corollaries~4.59 and 4.90]{krappdiss}, which encouraged me to compose this note, as well as the anonymous referee for valuable comments that improved the presentation of this note.
\end{acknowledgements}

	\begin{footnotesize}
		
	\end{footnotesize}
	
\end{document}